\documentclass[a4,12pt]{article}
\usepackage{amscd}
\usepackage{arydshln}
\usepackage{epsfig}
\usepackage{blkarray}
\usepackage{multirow}
\usepackage{amsxtra}
\usepackage{amssymb, amsmath}
\usepackage{amsthm}
\usepackage[mathscr]{eucal}
\usepackage{indentfirst}
\usepackage{amstext}
\usepackage{textcomp}
\newtheorem{theorem}{Theorem}[section]

\newtheorem{corollary}[theorem]{Corollary}

\begin{document}
\normalsize

\title{Infinite Symmetric Matrices over ${\mathbb Z}_2$ and the Lights Out Problem.}
\date{}

\author{Daniel Gon\c calves\footnote{Partially supported by CNPq, email:daemig@gmail.com} \\ \and Maria Inez Cardoso Gon\c calves\footnote{email:minezcg@gmail.com}}
\maketitle
\begin{center}
Departamento de Matematica - Universidade Federal de Santa Catarina - Trindade - Florian\'opolis - SC - 88.040-900 - Brazil.
\end{center}


\normalsize

\bigskip

\begin{abstract}

We show, using a hybrid analysis/linear algebra argument, that the diagonal vector of an infinite symmetric matrix over ${\mathbb Z}_2$ is contained in the range of the matrix. We apply this result to an extension, to the countable infinite case, of the Lights Out problem.

\end{abstract}

\medskip\noindent
{\it Key words}: Infinite matrices, lights out, symmetric, diagonal, range.

\smallskip\noindent  {\it AMS Subject Classifications:} 05C50, 15B33.

\section{Introduction:}

The Lights Out Problem arises in mathematics motivated by the electronic game Lights Out, which consists of 25 lighted buttons displayed as a $5 \times 5$ matrix. The game begins with some of the lights on and some of the lights off. At each step of the game, the player has to press one of the buttons, doing so changes the state of the pressed button, as well as the state of all the rectilinear adjacent buttons, that is, if the light is on then it is turned off and vice-versa. The goal of the game is to turn off all the lights. It is possible to develop a complete strategy for the game using linear algebra, more specifically using Gauss-Jordan elimination and the null and column space of a matrix. More details can be found in \cite{andersonfeil}.

The problem can also be seen as a $5\times 5$ graph, with a binary state assigned to each vertex, namely $0$ for light off and $1$ for light on. Changing the state of one vertex, changes the state of all the adjacent vertices. If we start the game with all the lights on, then it is possible to get the configuration where all the lights are off, by clicking on a subset of the vertices, see for example \cite{caro,minevich}. 

In \cite{minevich} Minevich proved that the range of a symmetric matrix over the field ${\mathbb Z}_2$ contains the vector of its diagonal elements and applied this result to a generalization of the Lights Out Problem (see \cite{andersonfeil, martinflores}), namely the case when  clicking some of the vertices does not actually affect the state of those vertices, but only the states of those adjacent to them.

The theorem in this note is a generalization of Minevich´s result to the case of an infinite symmetric matrix, with only a finite number of 1's in each line, over ${\mathbb Z}_2$. Such matrix could be, for example, the adjacency matrix of a countable infinite graph such that each vertex has finite degree (valency), and so this result can be used to solve an extension, to the infinite number of vertices case, of the Lights Out problem.

Before we proceed to the proof of the main theorem we recall the relevant results of \cite{minevich} below.

\section{Finite Symmetric Matrices over ${\mathbb Z}_2$}

The following result of \cite{minevich} will be used in the proof of theorem \ref{Oteorema}. 

\begin{theorem}\label{teo1} \cite{minevich} Let $A=[a_{ij}]$, $1\leq i \leq N$, $1\leq j \leq N$ be a symmetric matrix over ${\mathbb Z}_2$ and let $d=(a_{11},a_{22}, \cdots, a_{NN})^T$ be its diagonal vector. Then the range of $A$ contains the vector $d$.
\end{theorem}

This theorem is proved using the algebraic structure of ${\mathbb Z}_2$ and linear algebra results, more precisely the author uses the fact that the null space is the orthogonal complement of the range of a matrix.

\begin{corollary} \cite{minevich} \label{corominevich} Let $G$ be a finite undirected graph, where there can be at most one edge between any two vertices and a vertex can be connected to itself. Let $v_1,\cdots, v_n$  be the vertices connected to themselves. Assume that clicking a vertex makes the vertices adjacent to it (including itself if it is connected to itself) switch their state. Then, starting with all vertices at the state of 0, we can choose vertices to click so that as a result, precisely  $v_1,\cdots, v_n$ are 1 and the rest are 0.
\end{corollary}

For the Lights Out problem, the corollary ensures that starting with all the lights on, it is possible to turn off all the lights, by clicking in a finite number of vertices. It also ensures a solution for a generalization of the Lights Out problem to the case where not all the vertices affect themselves. 
For example: Suppose that we have a graph with vertices colored with blue and white lights, the vertices with blue lights affect themselves and also those around them, the vertices with white lights do not affect themselves, only the ones around them. Starting with all the lights off, it is possible to turn on all the blue lights and turn off all the white lights.

\section{Infinite Symmetric Matrices over ${\mathbb Z}_2$}

\begin{theorem}\label{Oteorema}
Let $A=[a_{ij}]$, $1\leq i < \infty$, $1\leq j < \infty$, be an infinite symmetric matrix over ${\mathbb Z}_2$, with a finite number of 1´s on each line, and $d=(a_{11},a_{22}, \cdots)^T$ be ts diagonal vector. Then the range of $A$ contains the vector $d$.
\end{theorem}
{\hskip -0.6cm{\bf Proof:  }}

In order to show that the range of $A$ contains the vector $d$, we will construct a vector $x=(x_1, x_2, \cdots )^T$ such that $Ax=d$. It will be key for this construction a certain block structure of the matrix $A$, which we will make precise below. 

Let $n$ be a positive integer and $k_0=0$. Consider the first $n$ rows of $A$. Since each row of $A$ contains only a finite number of 1's, we can find a positive integer $k_1 > 0$ such that $a_{i,j}=0$ for $ 1 \leq i \leq n$ and $n+k_1+1 \leq j < \infty$, and, since $A$ is symmetric, we also have that $a_{i,j}=0$ for $ n+k_1+1 \leq i < \infty $ and $1 \leq j \leq n$. 

Now, consider the first $n+k_1$ rows of $A$. Again, since each row of $A$ contains only a finite number of 1's, we can find a positive integer $k_2 > k_1$, such that  $a_{i,j}=0$ for $1 \leq i \leq n+k_1$ and $n+k_2+1 \leq j < \infty$, and $a_{i,j}=0$ for $ n+k_2 +1 \leq i < \infty $ and $1\leq j \leq n+k_1$.

Proceeding this way, we can find positive integers $k_3<k_4<k_5< \ldots$ such that, for $s=3,4,\ldots$, $a_{ij}=0$ for $1\leq i \leq n+k_{s-1}$ and $n+k_s+1\leq j < \infty$, and $a_{ij}=0$ for $n+k_s+1\leq i < \infty$ and $1\leq j \leq n+k_{s-1}$.

To get a finer picture of the above partition of the matrix $A$ into blocks we name the sub matrices of $A$. For simplicity, from now on we denote the element $a_{ii}$ of the diagonal by $d_i$. So, let $D_n$ be the square block containing the first $n$ entries of the diagonal of $A$, $D_{n+k_1}$ be the square block containing the diagonal elements $d_{n+1}, d_{n+2}, \ldots, d_{n+k_1}$ of $A$ and, for $j=2,3,\ldots$, let $D_{n+k_j}$ be the square block containing the diagonal elements $d_{n+k_{j-1}+1}, \ldots, d_{n+k_j}$ of $A$, that is, 
$$D_n=\left[\begin{array}{cccc}
d_1 &  * & \cdots & *  \\
* & d_2  & \cdots & *  \\ 
\vdots & \vdots  & \ddots & \vdots \\
* & * & \cdots  & d_n  \\
\end{array}\right]
, \quad D_{n+k_1}=\left[\begin{array}{cccc}
d_{n+1} &  * & \cdots & *  \\
* & d_{n+2}  & \cdots & *  \\ 
\vdots & \vdots  & \ddots & \vdots \\
* & * & \cdots  & d_{n+k_1}  \\
\end{array}\right]
, \ldots,$$ where $*$ denotes either 0 or 1.
We also call the sub matrices of $A$ that may or may not contain only zeros by $B_1, B_2, B_3, \ldots$ and the sub matrices that certainly contain only zeros by $0$. With this in mind we get the following picture of $A$:

\[
\begin{blockarray}{ccccccc}
 & \quad \quad \quad \quad n &   \quad \quad n+k_1 & \quad \quad n+k_2 & \quad \quad n+k_3 & \quad \quad n+k_4 &  \quad \quad \quad\quad \text{ }\\ 
 \begin{block}{c[c|c|c|c|c|c]}
       &   &         &   & & &\\  
       &D_n& B_1     &0  & 0 &0 & \cdots \\
  n    &   &         &   &  &   & \\ \cline{2-7}
       &   &         &   &  &   & \\  
       &B_1&D_{n+k_1}&B_2& 0 &0  & \cdots   \\  
 n+k_1 &   &         &   &   &   & \\ \cline{2-7}
       &   &         &   &   &   & \\  
       &0  &B_2      &D_{n+k_2}& B_3 & 0 & \cdots    \\  
 n+k_2 &   &         &   &   &  &   \\ \cline{2-7}
       &   &         &   &   &  & \\  
       &0  &0        &  B_3& D_{n+k_3} & B_4 &\cdots    \\  
 n+k_3 &   &         &   &   &   & \\\cline{2-7}
       &   &         &   &   &  & \\  
       &0  &0        &  0& B_4 & D_{n+k_4}  &\cdots    \\  
 n+k_4 &   &         &   &   &   & \\  \cline{2-7}
       &   &         &   &   &  & \\  
       &\vdots  &\vdots        &  \vdots& \vdots & \vdots &\ddots    \\  
       &   &         &   &   &   & \\   
 \end{block}
\end{blockarray}
\]

Now that we have the block structure of the matrix $A$ well understood, we can proceed to define the vector $x$ such that $Ax=d$. To do this we look at the finite square blocks of the matrix, obtain solutions for some block equations and use an analysis diagonal argument to define $x$. The details follow below.

Consider the submatrix $A_{{n+k_2}}=[a_{i,j}]$, where $1\leq i \leq n+k_2$ and $1\leq j \leq n+k_2$, that is,

\[
A_{{n+k_2}}=\left[\begin{array}{ccc:c:c}
          & & &         &   \\  
&      D_n& &B_1     &0   \\
 & &        &         &    \\ \hdashline
        & & &         &    \\  
&       B_1& &D_{n+k_1}&B_2   \\  
& &          &         &    \\ \hdashline
     & &    &         &    \\  
    &    0 &  &B_2      &D_{n+k_2}    \\  
         &  & &        &      
\end{array}\right]\]

Since $A_{{n+k_2}}$ is a $({n+k_2}) \times ({n+k_2})$ symmetric matrix  over ${\mathbb Z}_2$, from theorem \ref{teo1}, we have that the main diagonal of  $A_{{n+k_2}}$ belongs to the range  of $A_{{n+k_2}}$. Therefore, there exists $z^2=(z^2_1, z^2_2, \ldots, z^2_{n+k_2})^T$ such that, $\displaystyle{A_{{n+k_2}}z^2=(d_1,d_2,\dots,d_{n+k_2})^T}$. Let $c_i^{(j)}$ denote the $i$-th column of $A$ cut at the $j$-th row, that is, $c_i^{(j)}$ is the column vector formed by the first $j$ entries of the $i$-th column of $A$. Notice that $$\displaystyle{\sum_{i=1}^{n+k_2}} z^2_i\cdot c_i^{(n+k_2)}=(d_1,d_2,\dots,d_{n+k_2})^T,$$ and, in particular, since $a_{ij}=0$ for $1\leq i \leq n$ and $j> n+k_1$, we also have that $$\displaystyle{\sum_{i=1}^{n+k_1}} z^2_i\cdot c_i^{(n)}=(d_1,d_2,\dots,d_{n})^T$$

Now consider the submatrix $A_{{n+k_3}}=[a_{i,j}]$, where $1\leq i \leq n+k_3$ and $1\leq j \leq n+k_3$, that is:

\[
A_{{n+k_3}}=\left[\begin{array}{ccc:c:c:c}
          & & &         &  & \\  
&      D_n& &B_1     &0  & 0 \\
 & &        &         &  &  \\ \hdashline
        & & &         &   & \\  
&       B_1& &D_{n+k_1}&B_2& 0   \\  
& &          &         &    &\\ \hdashline
     & &    &         &   & \\  
    &    0 &  &B_2      &D_{n+k_2} & B_3    \\  
         &  & &        &  &  \\   \hdashline
     & &    &         &   & \\  
    &    0 &  &0      &B_3 & D_{n+k_3}     \\  
         &  & &        &  &   
\end{array}\right]\]

Again, since $A_{{n+k_3}}$ is a $({n+k_3}) \times ({n+k_3})$ symmetric matrix  over ${\mathbb Z}_2$, we have that the main diagonal of  $A_{{n+k_3}}$ belongs to the range  of $A_{{n+k_3}}$, so there exists $z^3=(z^3_1, z^3_2, \ldots, z^3_{n+k_3})^T$ such that $\displaystyle{A_{{n+k_3}}z^3=(d_1,d_2,\dots,d_{n+k_3})^T}$, and hence $$\displaystyle{\sum_{i=1}^{n+k_3}} z^3_i\cdot c_i^{(n+k_3)}=(d_1,d_2,\dots,d_{n+k_3})^T,$$ 
$$\displaystyle{\sum_{i=1}^{n+k_2}} z^3_i\cdot c_i^{(n+k_1)}=(d_1,d_2,\dots,d_{n+k_1})^T$$
and  
 $$\displaystyle{\sum_{i=1}^{n+k_1}} z^3_i\cdot c_i^{(n)}=(d_1,d_2,\dots,d_{n})^T.$$

Proceeding this way, we obtain a sequence of vectors $z^2,z^3,z^4,\cdots$,  such that for every positive integer $l \geq 2$ and $1\leq j \leq l-1$,  $$\displaystyle{\sum_{i=1}^{n+k_{l-j}}} z^l_i\cdot c^{(n+k_{(l-j)-1})}_i=(d_1, \cdots, d_{n+k_{(l-j)-1}})^T.$$ 

In particular, we have that 
$\displaystyle{\sum_{i=1}^{n+k_1}} z^l_i\cdot c^{(n)}_i=(d_1, \cdots, d_n)^T$ for all $l\geq 2$, and so the first $n+k_1$ entries of $z^l$ form a solution to the equation $\displaystyle{\sum_{i=1}^{n+k_{1}}} w_i\cdot c^{(n)}_i=(d_1, \cdots, d_n)^T$, $w_i \in {\mathbb Z}_2$. Since there is only a finite number of possible solutions for this last equation, there exists a subsequence $(z^{n_j})$ of $(z^n)$ such that $z^{n_j}_i=z_i^{n_{j+1}}$, for all $i=1,\cdots, n+k_1$ and $j\in {\mathbb N}$, that is, the first $n+k_1$ entries of vectors in the subsequence agree. These are the first $n+k_1$ entries of the solution vector $x$, that is, we let $x_i=z^{n_1}_i$, for $i=1,\cdots, n+k_1$.


Next we drop the first term in the subsequence above and pass to remaining subsequence. Now, for this sequence we have that  $\displaystyle{\sum_{i=1}^{n+k_{2}}} z^l_i\cdot c^{(n+k_1)}_i=(d_1, \cdots, d_{n+k_1})^T$, that is, the first $n+k_2$ entries of $z^l$ form a solution to the equation $\displaystyle{\sum_{i=1}^{n+k_{2}}} w_i\cdot c^{(n+k_1)}_i=(d_1, \cdots, d_{n+k_1})^T$, $w_i \in {\mathbb Z}_2$. Again, since there is only a finite number of possible solutions for this last equation, there exists a subsequence $(z^{n_j})$ of $(z^n)$ such that $z^{n_j}_i=x_i$, for all $1 \leq i \leq n+k_1$, $j\in {\mathbb N}$, and $z^{n_j}_i=z_i^{n_{j+1}}$ for all $n+k_1<i\leq  n+k_2$, $j\in {\mathbb N}$. With this in mind, we define the entries of the solution vector $x$  between $n+k_1$ and $n+k_2+1$ by $x_i=z^{n_1}_i$, for $n+k_1<i\leq n+k_2$.

 
Finally, proceeding inductively we define the entries of the solution vector $x$ between $n+k_m$ and $n+k_{m+1}+1$, for any $m\in {\mathbb N}$, and it is clear that the infinite vector $x$ defined this way is such that $Ax=d$ as desired.

\hfill{$\square$}\vskip 1pc

 Turning back to  the Lights Out problem,  we can use the above theorem to solve the following extension of the problem:  Let $\mathbf{G}$  be a countable infinite graph with finite degree vertices and such that there exists at most one edge between any two vertices and a vertex can be connected to itself. Is it possible, starting with all the lights out, to click on a (possibly) infinite number of vertices, and end up with only the vertices connected to themselves turned on?

In light of theorem \ref{Oteorema} one can easily adapt the proof given for the finite case in \cite{minevich} to the infinite case and so the question above is answered affirmatively. For completeness we make a precise statement below.

\begin{corollary}
Let $\mathbf{G}$ be a countable infinite graph such that each vertex has finite degree, there exists at most one edge between any two vertices and a vertex can be connected to itself. Let $(v_i)_{i\in \mathbb{N}}$ be an enumeration of the vertices of $\mathbf{G}$ and let  $(v_{i_1},v_{i_2},v_{i_3},\ldots)$ be the vertices connected to themselves. Then starting with all the lights off, it is possible, by clicking on a (possibly) infinite number of vertices, to turn on only the lights of $(v_{i_1},v_{i_2},v_{i_3}\ldots)$.
\end{corollary}


\end{document}